\documentclass[12pt]{article}

\usepackage{theorem,amssymb,amsmath}

\advance \topmargin by -\headheight
\advance \topmargin by -\headsep
\evensidemargin \oddsidemargin
\author{J.-P. Allouche\thanks{The author was partially supported by the ANR
project ``FAN'' (Fractals et Num\'eration), ANR-12-IS01-0002.} \\
CNRS, Institut de Math\'ematiques de Jussieu-PRG \\
\'Equipe Combinatoire et Optimisation \\
Universit\'e Pierre et Marie Curie, Case 247 \\
4 Place Jussieu \\
F-75252 Paris Cedex 05 France \\
{\tt jean-paul.allouche@imj-prg.fr}
\and
Thomas Baruchel \\
Lyc\'ee naval \\
Avenue de l'\'Ecole navale \\
F-29240 Brest Arm\'ees, France \\
{\tt baruchel@riseup.net}
}

\title{Variations on an error sum function for the convergents of some powers of $e$}

\date{ }

\def \proof{\bigbreak\noindent{\it Proof.\ \ }}

\def \endpf{{\ \ $\Box$ \medbreak}}
\DeclareMathOperator\erf{erf}

\newtheorem{theorem}{Theorem}

\newtheorem{corollary}{Corollary}

\theorembodyfont{\rm}

\newtheorem{remark}{Remark}
\newtheorem{example}{Example}

\begin{document}

\maketitle

\begin{abstract}
Several years ago the second author playing with different ``recognizers of real
constants'', e.g., the LLL algorithm, the Plouffe inverter, etc. found the following 
formula empirically. Let $p_n/q_n$ denote the $n$th convergent of the continued
fraction of the constant $e$. Then
$$
\sum_{n \geq 0} |q_n e - p_n| = 
\frac{e}{4} \left(- 1 + 10 \sum_{n \geq 0} \frac{(-1)^n}{(n+1)! (2n^2 + 7n + 3)}\right).
$$
The purpose of the present paper is to prove this formula and to give similar formulas
for some powers of $e$. 

\medskip


\medskip

\noindent
{\bf MSC Classes}: 11A55, 11J70, 11B83, 11B75.

\end{abstract}

\section{Introduction}

Playing with the convergents of $e$, the second author discovered several years ago the formula
\begin{equation}\label{several}
\sum_{n \geq 0} |q_n e - p_n| = 
\frac{e}{4} \left(- 1 + 10 \sum_{n \geq 0} \frac{(-1)^n}{(n+1)! (2n^2 + 7n + 3)}\right).
\end{equation}
While trying to prove the formula rigorously we began being interested in the following quantity.
If $\alpha$ is a positive real number, and if $p_n/q_n$ is the $n$th convergent of its
continued fraction, the quantity $|q_n \alpha - p_n|$ tends rapidly to zero. Thus the series
$\sum_{n \geq 0} |q_n \alpha - p_n|$ converges. This series measures in some sense the ``global
approximation'' of $\alpha$ by its convergents. We then learned from J. Shallit that the quantity 
$\sum_{n \geq 0} |q_n \alpha - p_n| = \sum_{n \geq 0} (-1)^n(q_n \alpha - p_n)$ was investigated 
in several papers \cite{Elsner2011, ElsnerStein2011, ElsnerStein2012, Elsner2014}, where the study 
of the quantity $\sum_{n \geq 0} (q_n \alpha - p_n)$ (first defined in \cite{RidleyPetruska}) can 
also be found. 

\medskip

It is natural to ask whether the sum of the series $\sum_{n \geq 0} |q_n \alpha - p_n|$
can be expressed in terms of $\alpha$ without explicitly using the convergents, in particular 
in the case where $\alpha$ has a ``nice'' continued fraction expansion, e.g., when $\alpha$ 
is quadratic or when $\alpha = e$. 

\section{Quadratic numbers}

The case of quadratic numbers was addressed in \cite{Elsner2011}
(also see \cite{Elsner2014}).

\begin{theorem}[Elsner]
Let $p_n/q_n$ be the $n$th convergent of the continued fraction of $\alpha$. Then 
the series $\sum_{n \geq 0} (q_n \alpha - p_n) x^n$ converges absolutely at least
for $|x| < \frac{1+\sqrt{5}}{2}$ and
$$
\sum_{n \geq 0} (q_n \alpha - p_n) x^n \in {\mathbb Q}[\alpha](x).
$$
In particular (taking $x=-1$), $\sum_{n \geq 0} |q_n \alpha - p_n|$ belongs to ${\mathbb Q}[\alpha]$.
\end{theorem}

\begin{example}[Elsner]

\ { } 

\begin{itemize}

\item $\sum_{n \geq 0} |q_n \sqrt{7} - p_n|  = \frac{7 + 5\sqrt{7}}{14}\cdot$

\item For any integer $n \geq 1$ we have
$\sum_{n \geq 0} |q_n (\frac{n + \sqrt{4+n^2}}{2}) - p_n| = \frac{1}{\frac{n + \sqrt{4+n^2}}{2} - 1}\cdot$

\item In particular $\sum_{n \geq 0} |q_n (\frac{1 + \sqrt{5}}{2}) - p_n| =  \frac{1 + \sqrt{5}}{2}\cdot$

\end{itemize}

\end{example}

\section{Powers of $e$}

Euler \cite{Euler} proved that the continued fraction expansion of $e$ is
$[2, 1, 2, 1, 1, 4, 1, 1, 6, 1, 1, 8,...]$ (sometimes replaced by the
not really regular expression $[1, 0, 1, 1, 2, 1, 1, 4, 1, 1, 6, 1, 1, 8,...]$).
After Euler, a large number of papers contained the computation of continued fraction expansions 
for some expressions containing $e$ (typically certain powers of $e$ possibly multiplied by some 
rational numbers, or numbers like $\frac{e^{2/k} - 1}{e^{2/k} + 1}$), see in particular 
\cite{Hermite, Lambert, Hurwitz, Lehmer18, Davis45, Olds70, MW70, vdP96, Cohn2006, Osler2006, Komatsu2007, 
McLaughlin2008, Komatsu2008, Komatsu2009-tok, Komatsu2009-roum, McCabe2009, Komatsu2012, Hetyei2014}.

\bigskip

The fundamental theorem we will use here is due to Komatsu \cite[Theorem 6, first part]{Komatsu2009-tok}.
Komatsu's theorem contains several previous results.

\begin{theorem}[Komatsu]\label{Kom1}
Let $\ell \geq 2$ and $s \geq 1$ be two integers. Let $p_n/q_n$ be the $n$th convergent of the 
continued fraction of 
$$
s e^{1/(\ell s)} = 
[s, \ell - 1, 1, 2s-1, 3\ell - 1, 1, 2s-1, 5\ell - 1, 1, 2s-1, \cdots, (2k-1 )\ell - 1 , 1, 2s-1, \cdots].
$$
Then for $n \geq 0$

\medskip

$p_{3n} - se^{1/(\ell s)}q_{3n} = 
- \displaystyle\frac{1}{(\ell s)^{n+1}} \int_0^1 \frac{x^n(x-1)^n}{n!} s e^{x/(\ell s)} \mbox{\rm d}x$

\bigskip

$p_{3n+1} - se^{1/(\ell s)}q_{3n+1} = 
\displaystyle\frac{1}{s(\ell s)^{n+1}} \int_0^1 \frac{(x+s-1)x^n(x-1)^n}{n!} s e^{x/(\ell s)} \mbox{\rm d}x$

\bigskip

$p_{3n+2} - se^{1/(\ell s)}q_{3n+2} = 
\displaystyle\frac{1}{s(\ell s)^{n+1}} \int_0^1 \frac{x^n(x-1)^{n+1}}{n!} s e^{x/(\ell s)} \mbox{\rm d}x$

\bigskip

\noindent
Let $s \geq 1$ be an integer. Let $p_n^*/q_n^*$ be the $n$th convergent of the continued fraction of 
$$
s e^{1/s} = 
[s+1, 2s-1, 2, 1, 2s-1, 4, 1, \cdots, 2s-1, 2k, 1, \cdots].
$$
Then $p_n^*/q_n^* = p_{n+2}/q_{n+2}$ with $p_n/q_n$ is as above. More precisely for $n \geq 0$

\medskip

$p_{3n}^* - se^{1/s}q_{3n}^* = 
\displaystyle\frac{1}{s^{n+2}} \int_0^1 \frac{x^n(x-1)^{n+1}}{n!} s e^{x/s} \mbox{\rm d}x$

\bigskip

$p_{3n+1}^* - se^{1/s}q_{3n+1}^* = 
- \displaystyle\frac{1}{s^{n+2}} \int_0^1 \frac{x^{n+1}(x-1)^{n+1}}{(n+1)!} s e^{x/s} \mbox{\rm d}x$

\bigskip

$p_{3n+2}^* - se^{1/s}q_{3n+2}^* = 
\displaystyle\frac{1}{s^{n+3}} \int_0^1 \frac{(x+s-1)x^{n+1}(x-1)^{n+1}}{(n+1)!} s e^{x/s} \mbox{\rm d}x$

\end{theorem}

\bigskip

Using Komatsu's result we can prove the following theorem. First recall that the ``error function''
$\erf$ is defined by $\displaystyle\erf(x) := \frac{2}{\sqrt{\pi}} \int_0^x e^{-t^2} \mbox{\rm d}t$.

\medskip

\noindent
{\bf Note} \ The name ``error sum function'' (or ``error-sum function'') that goes back at least
to \cite{RidleyPetruska} should not be confused with the name ``error function''. To (try to) avoid any
ambiguity, we will always write for the latter ``error function $\erf$''.

\begin{theorem}\label{AB1}
Let $\ell \geq 2$ and $s \geq 1$ be two integers. Let $p_n/q_n$ be the $n$th convergent of the
continued fraction of
$$
s e^{1/(\ell s)} = 
[s, \ell - 1, 1, 2s-1, 3\ell - 1, 1, 2s-1, 5\ell - 1, 1, 2s-1, \cdots, (2k-1 )\ell - 1 , 1, 2s-1, \cdots].
$$
Then
$$
\sum_{n \geq 0} |p_n - s e^{1/(\ell s)}q_n| = e^{1/\ell s} \sqrt{\frac{\pi s}{\ell}} \erf(1/\sqrt{\ell s}).
$$

\bigskip

\noindent
Let $s \geq 1$ be an integer. Let $p_n^*/q_n^*$ be the $n$th convergent of the continued fraction of
$$
s e^{1/s} = [s+1, 2s-1, 2, 1, 2s-1, 4, 1, \cdots, 2s-1, 2k, 1, \cdots].
$$
Then
$$
\sum_{n \geq 0} |p_n^* - s e^{1/s}q_n^*| = e^{1/s} \sqrt{\pi s} \erf(1/\sqrt{s})
+ s(1-e^{1/s}) - 1.
$$
\end{theorem}

\proof It suffices to use Komatsu's theorem (Theorem~\ref{Kom1} above) after writing
$$
\sum_{n \geq 0} |p_n - s e^{1/\ell s} q_n| =
\sum_{0 \leq j \leq 2} \sum_{n \geq 0} |p_{3n+j} - s e^{1/\ell s} q_{3n+j}|
$$
respectively
$$
\sum_{n \geq 0} |p_n^* - s e^{1/s} q_n^*| = |p_0^* - s e^{1/s} q_0| +
\sum_{1 \leq j \leq 3} \sum_{n \geq 0} |p_{3n+j}^* - s e^{1/s} q_{3n+j}^*|. \ \ \Box
$$

\bigskip

We deduce the following corollary.

\begin{corollary}\label{cor1}
Let $p_n^*/q_n^*$ be the $n$th convergent of the continued fraction of $e$
(recall that $e = [2, 1, 2, 1, 1, 4, 1, \cdots, 1, 2n, 1, \cdots]$). Then
$$
\sum_{n \geq 0} |p_n^* - e q_n^*| = 2 e \int_0^1 e^{-t^2} dt - e 
= e \sqrt{\pi} \erf(1) - e.
$$

\medskip

\noindent
Let $p_n/q_n$ be the $n$th convergent of the continued fraction of $e^{1/\ell}$ (with $\ell \geq 2$).
Then 
$$
\sum_{n \geq 0} |p_n - e^{1/\ell} q_n| = e^{1/\ell} \sqrt{\frac{\pi}{\ell}} \erf(1/\sqrt{\ell}).
$$
\end{corollary}

\begin{remark}
The first result in Corollary~\ref{cor1} above was already obtained by Elsner in
\cite[p.\ 2]{Elsner2011}.
\end{remark}

\bigskip

Now we prove Formula~(\ref{several}).

\begin{corollary}\label{cor-for}
Let $A(\ell, s)$ be defined for positive reals $\ell$ and $s$ by
$$
A(\ell, s) := \sum_{n \geq 0} \frac{(-1)^n}{(n+1)! (2n^2 + 7n + 3) (\ell s)^n}\cdot
$$
Then 
$$
A(\ell, s) = - \frac{3}{10} \ell s + \frac{1}{5} \ell s (2 - \ell s - \ell^2 s^2) e^{-1/\ell s} +
\frac{4}{5} \int_0^1 e^{-t^2/\ell s} \mbox{\rm d}t
$$
In particular 
$$
A(1,1) = -\frac{3}{10} + \frac{4}{5} \int_0^1 e^{-t^2} \mbox{\rm d}t.
$$
so that
$$
\sum_{n \geq 0} |q_n e - p_n| = 2 \int_0^1 e^{-t^2} \mbox{\rm d}t - e
= \frac{e}{4} (- 1 + 10 A(1,1))
$$
\end{corollary}

\proof The proof is easy. First write 
$$\frac{1}{2n^2 + 7n + 3} = \frac{2}{5(2n+1)} - \frac{1}{5(n+3)}\cdot
$$
Then introduce the series
$$
\sum_{n \geq 0} \frac{(-1)^{n+1} x^{2n+1}}{(n+1)!(\ell s)^{n+1} (2n+1)}
\ \ \mbox{\rm and} \ \
\sum_{n \geq 0} \frac{(-1)^{n+1} x^{n+3}}{(n+1)!(\ell s)^{n+1} (n+3)}\cdot
$$
The derivative of these series are easily computed. We then need their values at $x=1$. \endpf

\begin{remark}
It is immediate to use the values of $A(\ell, s)$ to obtain similar formulas for $se^{\ell s}$. 
We also note we first thought that the quantity $(2n^2 + 7n + 3)$ was somehow crucial in 
Formula~(\ref{several}): there might ever have been (though it would have been quite surprising) 
a link with the number of independent parameters of the orthosymplectic group ${\rm OSP}(3,2n)$ 
which is precisely $(2n^2 + 7n + 3)$ (see, e.g., \cite[p.\ 223]{Chamseddine}). 
But this quantity is not crucial; compare with Formula~(\ref{easier}) given below which can be proved 
by using a step of the proof of Corollary~\ref{cor-for} above with $s=\ell=1$ and $x=1$:
$$
\sum_{n \geq 0} \frac{(-1)^{n+1}}{(n+1)!(2n+1)}
= 1 - e^{-1} - 2\int_0^1 e^{-t^2} \mbox{\rm d}t.
$$
This implies
\begin{equation}\label{easier}
\sum_{n \geq 0} |p_n^* - e q_n^*| = e \sum_{n \geq 0} \frac{(-1)^n}{(n+1)!(2n+1)} - 1.
\end{equation}

\end{remark}

\begin{remark}
The value of $\sum_{n \geq 0} |p_n - \alpha q_n]$ for $\alpha$ equal to one of the above real 
numbers can also be expressed as another kind of series. Namely a classical series for the error 
function $\erf$ (see, e.g., \cite[7.1.6, p.\ 297]{Abr-Ste} reads
$$
\erf(z) = \frac{2}{\sqrt{\pi}}e^{-z^2} 
\sum_{n \geq 0} \frac{2^n}{1 \times 3 \times 5 \cdots \times (2n+1)}z^{2n+1}.
$$
Using Corollary~\ref{cor1} and the notation therein, this gives in particular the following formulas
\begin{equation}\label{other-e}
\sum_{n \geq 0} |p_n^* - e q_n^*| 
= \sum_{n \geq 0} \frac{2^{n+1}}{1 \times 3 \times 5 \cdots \times (2n+1)} - e
= \sum_{n \geq 0} \frac{2^{2n+1} \, n!}{(2n+1)!} - e
\end{equation}
and, for any integer $\ell \geq 2$,
\begin{equation}\label{other-e-powers}
\sum_{n \geq 0} |p_n - e^{1/\ell} q_n| = 
\sum_{n \geq 0} \frac{2^{n+1}}{\ell^{n+1} \, (1 \times 3 \times 5 \cdots \times (2n+1))}
= \sum_{n \geq 0} \frac{2^{2n+1} \, n!}{\ell^{n+1} \, (2n+1)!}\cdot
\end{equation}
Note that the second author obtained Formula~(\ref{other-e-powers}) empirically.
Also note that the digits of the decimal expansion of the right side of Equation~(\ref{other-e})
(up to the $-e$ term) is given in \cite{oeis} as A125961, and that the expansions of the right 
side of Equation~(\ref{other-e-powers}) above for $\ell = 2$ and $\ell = 4$ are given in \cite{oeis} 
as A060196 and A214869 respectively.
\end{remark}

\section{More fun with the error sum function}

Formulas similar to the formulas in the previous section can be stated by using results on the 
convergents for continued fractions with ``regular'' patterns, in particular at least for 
(some of) the so-called Hurwitz continued fractions, sometimes also called (regular) continued 
fractions of Hurwitzian type, see, e.g., \cite{McLaughlin2008}.
We simply list below results that can be used to yield nice formulas for the error sum function 
we considered. They give in terms of integrals for some reals $\alpha$ and their convergents $p_n/q_n$ 
the quantity $p_{an+b} - \alpha q_{an+b}$ (for any $b$ in a complete system of residues modulo $a$)
and they are due to Komatsu.

\begin{itemize}

\item for $\alpha = \frac{e^{1/\ell s}}{s}$, with $s$ and $\ell$ any two integers $\geq 2$,
and $\alpha = \frac{e^{1/s}}{s}$, with $s \geq 2$, integral expressions for 
$p_{3n+j} - \alpha q_{3n+j}$ with $j \in \{0,1,2\}$ are given in 
\cite[Theorem 3, second part]{Komatsu2009-tok}.

\item for $\alpha = e^{2/s}$, with $s \geq 3$ and odd, integral expressions for 
$p_{5n+j} - \alpha q_{5n+j}$ with $j \in \{0, 1, 2, 3, 4\}$ are given in \cite{Komatsu2007};

\item for $\alpha = \frac{e^{1/(3s+1)}}{3}$ (resp.\ $\alpha = \frac{e^{1/(3s+2)}}{3}$), integral 
expressions for $p_{9n+j} -\alpha q_{9n+j}$ with $j \in \{-6, -5, -4, -3, -2, -1, 0, 1, 2\}$ are 
given in \cite{Komatsu2009-tok};

\item for $\alpha = \sqrt{\frac{v}{u}} \tanh \frac{1}{\sqrt{uv}}$, integral expressions for
$p_{2n-1} - \alpha q_{2n-1}$ and $p_{2n} - \alpha q_{2n}$ are given in \cite{Komatsu2009-roum};

\item for $\alpha = \sqrt{\frac{v}{u}} \tan \frac{1}{\sqrt{uv}}$, integral expressions for
$p_{4n-j} - \alpha q_{4n-j}$ with $j \in \{0, 1, 2, 3\}$ are given in \cite{Komatsu2009-roum}.

\end{itemize}

\bigskip

The last result we would like to cite here is a nice particular case of a theorem
of Hetyei \cite[Theorem 2.9]{Hetyei2014} (also see \cite[p.\ 21]{Hetyei2014})
which could be used to compute the error sum function for 
$\alpha = \frac{4(11\sin(1/2) - 6\cos(1/2))}{53\cos(1/2)-97\sin(1/2)}\cdot$

\begin{theorem}[Hetyei]
We have the following continued fraction expansion
$$
\frac{4(11\sin(1/2) - 6\cos(1/2))}{53\cos(1/2)-97\sin(1/2)} = [4,3,4,4,4,5,4,6,4,7,4,...].
$$
\end{theorem}

\section{Conclusion} The error sum function of some other continued fractions with 
``regular'' patterns could probably be studied. Another appealing possibility is the 
definition and study of error sum functions similar continued fractions in the function 
field case (see in particular \cite{Thakur1, Thakur2, Thakur3}). Finally we give a last
relation that the second author discovered empirically: we did not locate it in the literature
(yet) and did not prove it (yet)
$$
\int_0^1 e^{-t^2} \mbox{\rm d}t = 3/8 + 
\cfrac{5/4}{3+\cfrac{9}{21 + \cfrac{288}{63 + \cfrac{\ddots}{\ddots + \cfrac{n(n+2)^2(2n-1)^2}{
(2n+5)(n^2+n+1) + \ddots}}}}}
$$

\bigskip

\noindent
{\bf Acknowledgements} We want to warmly thank J. Shallit who indicated to us
the papers \cite{Elsner2011, ElsnerStein2011, ElsnerStein2012, Elsner2014}.

\bigskip

\noindent
{\bf Addendum} After we posted the first version of this paper on ArXiv, C. Elsner kindly 
sent us the preprint \cite{ElsnerKlauke}, which was written a few months ago, and where 
the reader can find some results in another direction (limit formulas, differential equations, 
algebraic independence results, relations to Hall's theorem) but also a proof of the result 
in Corollary~\ref{cor1}
$$
\forall \ell \geq 2, \ \sum_{n \geq 0} |p_n - e^{1/\ell} q_n| = e^{1/\ell} \sqrt{\frac{\pi}{\ell}} \erf(1/\sqrt{\ell}).
$$


\begin{thebibliography}{99}

\bibitem{Abr-Ste} M. Abramowitz, I. A. Stegun, Handbook of Mathematical Functions With Formulas, 
Graphs, and Mathematical Tables, National Bureau of Standards Applied Mathematics Series {\bf 55},
10th printing, 1972.

\bibitem{Chamseddine} A. H. Chamseddine, Massive supergravity from spontaneously breaking 
orthosymplectic gauge symmetry, {\it Ann. Physics\,} {\bf 113} (1978) 219--234.

\bibitem{Cohn2006} H. Cohn, A short proof of the simple continued fraction expansion of $e$, 
{\it Amer. Math. Monthly\,} {\bf 113} (2006) 57--62. 

\bibitem{Davis45} C. S. Davis, On some simple continued fractions connected with $e$, 
{\it J. London Math. Soc.} {\bf 20} (1945) 194--198. 

\bibitem{Elsner2011} C. Elsner, Series of error terms for rational approximations of irrational 
numbers, {\it J. Integer Seq.} {\bf 14} (2011), Article 11.1.4, 20 p. 

\bibitem{Elsner2014} C. Elsner, On error sums of square roots of positive integers with applications 
to Lucas and Pell numbers, {\it J. Integer Seq.} {\bf 17} (2014), Article 14.4.4, 21 p. 

\bibitem{ElsnerKlauke} C. Elsner, A. Klauke, Errorsums for the values of the exponential function,
{\it Bericht Nr. 02014/01}, 1--19;  Technische Informationsbibliothek Hannover: RS 8153 (2014,1).

\bibitem{ElsnerStein2011} C. Elsner, M. Stein, On error sum functions formed by convergents of real 
numbers, {\it J. Integer Seq.} {\bf 14} (2011), Article 11.8.6, 14 p. 

\bibitem{ElsnerStein2012} C. Elsner, M. Stein, On the value distribution of error sums for approximations 
with rational numbers, {\it Integers\,} {\bf 12} (2012), Article A66, 28 p. 

\bibitem{Euler} L. Euler, De fractionibus continuis dissertatio, {\it Comm. Acad. Sci. Petropol.}
{\bf 9} (1744) 98--137. Available at 
{\tt http:/$\!$/eulerarchive.maa.org/docs/originals/E071.pdf} or at
{\tt http:/$\!$/www.math.dartmouth.edu/$\sim$euler/docs/originals/E071.pdf}

\noindent
English translation: M. F. Wyman, B. F. Wyman, An essay on continued fractions, {\it Math. Systems 
Theory\,} {\bf 18} (1985) 295--328.

\bibitem{Hermite} C. Hermite, Sur la fonction exponentielle, {\it C. R. Acad. Sci.} {\bf 77} 
(1873) 18--24, 74--79, 226--233, and 285--293; also in {\OE}uvres de Charles Hermite, publi\'ees
sous les auspices de l'Acad\'emie des sciences par \'Emile Picard, vol. 3, Gauthier-Villars, 
Paris, 1912, pp. 150--181. Available at 
{\tt https:/$\!$/archive.org/details/oeuvresdecharles03hermuoft} 

\bibitem{Hetyei2014} G. Hetyei, Hurwitzian continued fractions containing a repeated constant and 
an arithmetic progression, {\it SIAM J. Discrete Math.} {\bf 28} (2014) 962--985. 

\bibitem{Hurwitz} A. Hurwitz, \"Uber die Kettenbr\"uche, deren Teilnenner arithmetische Reihen bilden,
{\it Z\"urich. Naturf. Ges.} {\bf 41} 2nd Part (1896) 34--64.

\bibitem{Komatsu2007} T. Komatsu, A proof of the continued fraction expansion of $e^{2/s}$,
{\it Integers} {\bf 7} (2007), Article A30, 8 p.

\bibitem{Komatsu2008} T. Komatsu, More on Hurwitz and Tasoev continued fractions,
{\it Sarajevo J. Math.} {\bf 4} (2008) 155--180. 

\bibitem{Komatsu2009-roum} T. Komatsu, Diophantine approximations of $\tanh$, $\tan$ and linear 
forms of $e$ in terms of integrals, {\it Rev. Roum. Math. Pures Appl.} {\bf 54} (2009) 223--242.

\bibitem{Komatsu2009-tok} T. Komatsu, A diophantine approximation of $e^{1/s}$ in terms of integrals, 
{\it Tokyo J. Math.} {\bf 32} (2009) 159--176.

\bibitem{Komatsu2012} T. Komatsu, Some exact algebraic expressions for the tails of Tasoev continued 
fractions, {\it J. Aust. Math. Soc.} {\bf 92} (2012) 179--193.

\bibitem{Lambert} J. H. Lambert, M\'emoires sur quelques propri\'et\'es remarquables des quantit\'es 
transcendantes, circulaires et logarithmiques, {\it M\'emoires de l'Acad\'emie royale des sciences 
de Berlin} {\bf 17} (1761/1768) 265--322. Available at 
{\tt http:/$\!$/www.kuttaka.org/$\sim$JHL/L1768b.pdf}

\bibitem{Lehmer18} D. N. Lehmer, Arithmetical theory of certain Hurwitzian continued fractions, 
{\it Amer. J. Math.} {\bf 40} (1918) 375--390.
 
\bibitem{McCabe2009} J. H. McCabe, On the Pad\'e table for $e^x$ and the simple continued fractions 
for $e$ and $e^{L/M}$, {\it Ramanujan J.} {\bf 19} (2009) 95--105. 

\bibitem{McLaughlin2008} J. Mc Laughlin, Some new families of Tasoevian and Hurwitzian continued 
fractions, {\it Acta Arith.} {\bf 135} (2008) 247--268. 

\bibitem{MW70} K. R. Matthews, R. F. C. Walters, Some properties of the continued fraction expansion of 
$(m/n)e^{1/q}$, {\it Proc. Cambridge Philos. Soc.} {\bf 67} (1970) 67--74. 

\bibitem{Olds70} C. D. Olds, The simple continued fraction expansion of $e$, {\it Amer. Math. 
Monthly\,} {\bf 77} (1970) 968--974.

\bibitem{oeis} On-Line Encyclopedia of Integer Sequences, available electronically
at {\tt http:/$\!$/oeis.org}

\bibitem{Osler2006} T. J. Osler, A proof of the continued fraction expansion of $e^{1/M}$, 
{\it Amer. Math. Monthly\,} {\bf 113} (2006) 62--66. 

\bibitem{vdP96} A. J. van der Poorten, Continued fraction expansions of values of the exponential 
function and related fun with continued fractions, {\it Nieuw Arch. Wisk.} {\bf 14} (1996) 221--230. 

\bibitem{RidleyPetruska} J. N. Ridley, G. Petruska, The error-sum function of continued fractions,
{\it Indag. Math.} {\bf 11} (2000) 273--282.

\bibitem{Thakur1} D. S. Thakur, Continued fraction for the exponential for ${\bf F_q[T]}$,
{\it J. Number Theory\,} {\bf 41} (1992) 150--155.

\bibitem{Thakur2} D. S. Thakur, Exponential and continued fractions, {\it J. Number Theory\,}
{\bf 59} (1996) 248--261. 

\bibitem{Thakur3} D. S. Thakur, Patterns of continued fractions for the analogues of $e$ and related 
numbers in the function field case, {\it J. Number Theory\,} {\bf 66} (1997) 129--147. 

\end{thebibliography}
\end{document}